\documentclass[a4paper,fleqn]{cas-sc}

\usepackage[numbers]{natbib}

\def\tsc#1{\csdef{#1}{\textsc{\lowercase{#1}}\xspace}}
\tsc{WGM}
\tsc{QE}

\newtheorem{theorem}{Theorem}
\newtheorem{lemma}[theorem]{Lemma}
\newdefinition{rmk}{Remark}
\newproof{pf}{Proof}
\newproof{pot2}{Proof of Theorem \ref{thm2}}
\newproof{pot3}{Proof of Theorem \ref{thm3}}

\newtheorem{proposition}[theorem]{Proposition}
\newtheorem{conclusion}[theorem]{Corollary} 

\usepackage{tikz}
\usetikzlibrary{positioning}
\tikzset{main node/.style={circle,fill=blue!20,draw,minimum size=0.5cm,inner sep=0pt}}

\def\B{\mathcal B}
\def\Tau{\mathcal{T}}
\def\M{{\,\tt M\,}}
\def\FM{{\,\tt FM\,}}
\def\MR#1{MR \href{http://www.ams.org/mathscinet-getitem?mr=#1}{#1}}
\newcommand{\doi}[1]{\href{http://dx.doi.org/#1}{\texttt{doi:#1}}}
\DeclareMathOperator{\Gau}{Gau}
\DeclareMathOperator{\rank}{rank}
\DeclareMathOperator{\sign}{sign} 

\begin{document}
\let\WriteBookmarks\relax
\def\floatpagepagefraction{1}
\def\textpagefraction{.001}

\shorttitle{The $\alpha$-representation 
for the Tait coloring and for the $\chi_M(q)$}

\shortauthors{Eduard Lerner} 

\title [mode = title]{The $\alpha$-representation 
for the Tait coloring and for the
characteristic polynomial
of matroid}

\tnotemark[1] 

\tnotetext[1]{The research was supported by RSF (project No. 24-21-00158).}  

%
\author{Eduard~Lerner}[type=editor,
                        auid=000,bioid=1,
                        orcid=0000-0001-5129-5820
]



\begin{abstract}
Consider a finite field $\mathbb F_q$, $q=p^d$, where $p$ is an odd number.
Let $M=(E,r)$ be a regular matroid; denote by ${\mathcal B}$ the family of its bases, $\bar s(M;\alpha)=\sum_{B\in{\mathcal B}}\prod_{e\not\in B} \alpha_e$, where ${\alpha_e\in \mathbb F_q}$, $\alpha_e\neq 0$. Let a subset $A\equiv A(\alpha)$ in $E$ have the maximal cardinality and satisfy the condition $\bar s(M|A;\alpha)\neq 0$, while $r^*(\alpha)=|A|-r(E)$. Let us represent the value of the characteristic polynomial of the matroid $M$ at the point $q$ as the linear combination of Legendre symbols with respect to $\bar s(M|A;\alpha)$, whose coefficients are modulo equal to $1/q^{r^*(\alpha)/2}$. This representation generalizes the formula for a flow polynomial of a graph which was obtained by us earlier. The latter formula is an analog of the so-called $\alpha$-representation of vacuum Feynman amplitudes in the case of a finite field, which has inspired the Kontsevich conjecture (1997). The $\alpha$-representation technique is also applicable for expressing the number of Tait colorings for a cubic biconnected planar graph in terms of principal minors of the matrix of faces of this graph.
\end{abstract}

    
\maketitle

\section{Introduction}
\subsection{The historical background and the goal of the paper}
Since the very appearance of the theory of matroids, its significant problems have been related to finite fields;
it suffices to mention the critical problem stated by H.~Crapo and G.-C.~Rota~\cite{crapo}.
Nevertheless, in calculations of values of characteristic polynomials of matroids, one rarely immediately uses theoretical numerical ingredients. 

The idea of this paper has a long history. Let the symbol $\Tau(G)$ stand for the set of spanning trees of the graph~$G$. Consider sums
\begin{equation} 
\label{eq:2}
s(G;\alpha)=\sum_{T \in \Tau(G)}
  {\prod_{e\in E(T)} {\alpha_e } },
\end{equation}
where $\alpha_e$ are elements of the finite field ${\mathbb F}_q$.
In December 1997, when giving a talk at the Gelfand seminar at Rutgers University, Maxim Kontsevich proposed a conjecture that the number of nonzero values of~(\ref{eq:2}) for $\alpha  \in \mathbb{F}_q^E$ is a polynomial with respect to~$q$.
This conjecture was inspired by studying analogous sums (with complex-valued $\alpha_e$) in the quantum field theory. Though the conjecture was never published, it has aroused the interest of experts in combinatorics (see~\cite{stanleyArticle,chung,Stembridge}).
This conjecture allows a natural generalization for matroids, while all possible bases of a matroid represent an analog of $\Tau(G)$.
However, the increased interest to the Kontsevich conjecture and its generalization has negatively affected the development of such techniques, when sometime later this conjecture was refuted~\cite{belk}.
Nevertheless, accurately adjusting $\alpha$-representation techniques used in the quantum theory for fields of real values for the case of finite fields\footnote{See Remark~\ref{rem:alpha} (subsection~\ref{bog}) for the brief history of the term ``$\alpha$-representation''.}, we get a new representation for the flow polynomial of the graph~\cite{EJC}.

The goal of this paper is to generalize results obtained in the paper~\cite{EJC} for the case of arbitrary regular matroids and to apply the $\alpha$-representation technique for deriving explicit formulas in the case, which plays an important role in the graph theory.
Namely, we, in particular, consider the number of Tait colorings for a cubic biconnected planar graph. The fact that the number of Tait colorings for any such graph differs from zero is equivalent to the assertion of the Four Colors Theorem.

In obtaining the $\alpha$-representation for the characteristic polynomial of an arbitrary regular matroid the fact that the representation matrix is unimodular plays an important role.
In a non-regular case, we get a more complicated formula which depends on the representation matrix of the matroid over the field~$\mathbb F_q$.

The expression for the number of Tait colorings for a cubic planar graph~$G$ has appeared to be rather simple.
It is related to the matrix of faces of this graph, which (as far as we know) was not studied earlier. 
In this matrix, rows and columns correspond to all possible faces $F_i$ of the graph $G$. The $i,j$-th element of this matrix equals the sum of values of variables $\alpha(v)$, where $v$ is the vertex, which belongs to both faces $F_i$ and $F_j$, 
while $\alpha(v)$ are nonzero variables in the field $\mathbb F_3$ associated with the graph vertices 
(we treat these variables, which take on values $\pm 1$, as spins).
The number of Tait colorings depends on the sum of Legendre symbols of principal nonzero minors of this matrix with the maximal possible rank with respect to all possible $\alpha$.

The rest part of the paper has the following structure.
In the next subsection we state the main results, below we give an illustrative example for them.
Section~\ref{sect2} is devoted to auxiliary assertions used in the proof of the main results.
In particular, we state (and later prove) a more general theorem, which is valid for all matroids representable over the field~$\mathbb F_q$.
An important role in this research is played by (adduced in subsection~\ref{subcec2-5}) expressions for the characteristic function
of the dual matroid in terms of representation matrices of the initial matroid and properties of the generalized
Laplace--Kirchhoff matrix of an arbitrary matroid.
In subsection~\ref{subcec2-6}, we also recall the Heawood theorem for the number of Tait colorings; below we make use of this theorem.
In the application of the $\alpha$-representation in the case of a real field, one uses explicit formulas for calculating Gaussian integrals.
In the case of a finite field, in subsection~\ref{subcec2-7} we consider multidimensional Gaussian sums.
Section~\ref{sec3} is immediately devoted to the $\alpha$-representation. We first prove the theorem on the number of Tait colorings. Those, who are interested only in understanding this result, may restrict themselves to reading only corresponding parts of subsections~\ref{subsec1-2} and~\ref{subsec1-3} and, after getting acquainted with subsections~\ref{subcec2-6},~\ref{subcec2-7}, read this proof in Section~\ref{sec3}. 
We use it for clarifying a rather simple principle of the proof of the main result about the alpha-representation of matroids.
We believe that after a thorough preparation a reader would easily understand this idea. 
In Conclusion, we summarize the obtained results and discuss the further development of this work.    

\subsection{Denotations and statements of main results}
\label{subsec1-2}
Let us exactly state the main results of this paper.
We understand {\it a matroid} $M$ as a pair $(E,r)$, where $r\equiv r_M$ is the \textit{rank function} of $M$,
and $E$ is the \textit{ground set} of~$M$. 
Recall~\cite{Oxley} that $r$ 
is defined on all subsets $2^E$ of $E$, takes on nonnegative integer values in $\mathbb N_0$,
and possesses the following properties:\\
(1) If $A\subseteq E$, then $0 \leq r(A) \leq |A|$.\\
(2) If $B \subseteq A \subseteq  E$, then $r(B) \leq r(A)$.\\
(3) If $A$ and $B$ are subsets of $E$, then $r(A \cup B ) + r(A \cap B ) \leq r(A) + r(B )$.

We denote a singleton $\{e\}\subseteq E$ as $e$; if $r(e)=0$, then $e$ is said to be a {\it loop}. 
By definition, $r(M)\equiv r(E)$; the set $B$, $B\subseteq E$, is called a {\it base}, if $r(B)=|B|=r(M)$. Denote the totality of all bases of a matroid as~$\B\equiv \B_M$.

Let $A \subseteq E$. The restriction of $M = (E, r)$ to $A$ is the matroid $M |A$ with the ground set $A$ and the rank function coinciding with $r$ on $2^A$.

Recall that the characteristic polynomial $\chi_M$ of the matroid~$M$ obeys the formula
$$
\chi_M(x)=\sum_{A:A\subseteq E} (-1)^{|A|} x^{r(E)-r(A)}.
$$
As one can easily see, if $M$ contains loops, then $\chi_M(x)\equiv 0$.

A matroid is said to be {\it representable} over a field $\mathbb F$, if there exists a matrix $\M$ such that its elements belong to the field~$\mathbb F$, numbers of its columns are elements of the set~$E$, and $r(A)$ equals the rank of the set of columns of the matrix~$\M$ with indices in~$A$.
A matroid is said to be {\it regular}, if it is representable over any field $\mathbb F$ (see, for example, \cite{Oxley,reiner1} for more details).

We consider {\it the field}~$\mathbb F_q$ of an odd characteristic, i.e., $q=p^d$, where $p$ is an odd prime, and $d\in\mathbb N$. Denote by $\eta$ {\it the multiplicative quadratic character} of the field~$\mathbb{F}_q$: $\eta(0)=0$, in other cases, $\eta(x)=1$ or $\eta(x)=-1$ depending on whether $x$ is a square in the field $\mathbb F_q$ or not.
For $d=1$ we have $\eta(x)=\left(\frac{x}{p}\right)$, where $\left(\frac{x}{p}\right)$ is {\it the Legendre symbol} of the residue field modulo prime~$p$. 
Let us define a function $g(q,m)$, where $q$ is the number of elements in the mentioned field, $m\in\mathbb N$, by the formula
\[
 g(q,m) = \left\{ \begin{array}{ll}
   1/q^{m/2},& \textrm{if}~p \bmod 4= 1,\,  m\bmod 2= 0,  \\
   1/(-q)^{m/2},&  \textrm{if}~p \bmod 4=3,\,  m\bmod 2= 0,\\
   0,& \textrm{if}~ m\bmod 2=1.
     \end{array} \right.
\]

Let us associate each element of the matroid $e$, $e\in E$, with a nonzero value $\alpha_e$ in the field~$\mathbb F_q$ and assume that
$$
s(M)\equiv s(M;\alpha)=\sum_{B\in\B}\prod_{e\in B} \alpha_e,\qquad \bar s(M) \equiv\bar s(M;\alpha)=\sum_{B\in\B}\prod_{e\not \in B} \alpha_e.
$$
Here the product over the empty set equals one. 
In subsection~\ref{bog}, we discuss the historical background of the choice of the denotation $\alpha$. Let us put 
$$
A^*(M;\alpha)\equiv A^*=\mathop{\rm argmax}\limits_A \{|A|: s(M|A)\neq 0 \},\quad   r^*(M;\alpha)=|A^*|-r(M).
$$

Note that with nonzero $\alpha_e$ the condition $s(M|A;\alpha)\neq 0$ is equivalent to the inequality $\bar s(M|A;\alpha)\neq 0$.
The set $A^*(\alpha)$ is not defined uniquely, but the value of the function $\eta$ in the next theorem is independent of its choice.

\begin{theorem}
\label{thm1}
Put $q=p^d$, where $p$ is an odd prime number.
The following formula is valid for any regular matroid without loops:
\begin{equation} \label{eq:main1}
\chi_{M} (q) =
\sum
g(q,r^*(M;\alpha))\ 
\eta(\bar s(M|A^*));
\end{equation}
here the sum is calculated over everywhere nonzero values $\alpha$ in $\mathbb F_q^E$.
\end{theorem}

Note that one can seek for the subset $A^*$ through a narrower set, see Co\-rol\-la\-ry~\ref{conc:1} in subection~\ref{A*} for more details.

Let us now consider the result of another application of the same technique. Given a {\it simple biconnected planar cubic} graph $G=(V,E)$ (in this case, the biconnectivity of the graph is equivalent to the absence of bridges in it), we assume that the number of edges in the graph equals $3n$, while the number of vertices and faces is $2n$ and $n+2$, correspondingly, $n=2,3,\ldots$.
{\it The Tait coloring} is a coloring of edges in~$E$ in 3 colors such that all edges with a common vertex are colored differently.
The existence of such coloring for any graph~$G$ in the class under consideration is equivalent to the assertion of the Four Colors Theorem. We denote {\it the number of various Tait colorings} for the graph $G$ by the symbol $\chi'_G(3)$.

Denote by $\alpha(v)$ variable ({\it spins}) of vertices $v$ ($v\in V(G)$) which take on values $\pm 1$ in the field $\mathbb F_3$. 
Let $F_1,\ldots,F_{n+2}$ be all faces of the graph $G$. We understand the matrix of faces of the graph $G$ as an $(n+2)\times (n+2)$-matrix $\FM(\alpha)$, whose $i,j$-th element equals the sum of values $\alpha(v)$, where $v$ belongs to both faces $F_i$ and $F_j$. 

Evidently, the matrix $\FM(\alpha)$ is symmetric and degenerate (since the graph is cubic, the sum of elements in any row of this matrix equals zero in the field $\mathbb F_3$). Denote the rank of the matrix $\FM(\alpha)$ by the symbol $r(\FM(\alpha))$. 
Let $\FM'(\alpha)$ be an arbitrary principal minor $\FM(\alpha)$ of the order $\rank (\FM(\alpha))$ such that $\det  \FM'(\alpha)$ differs from zero. 
Let $\sign ( \FM(\alpha))=+1$, if $\det  \FM'(\alpha)=1$, and $\sign ( \FM(\alpha))=-1$ otherwise. 
In other words, $\sign ( \FM(\alpha))=\left(\frac{\det  \FM'(\alpha)}{3}\right)$.

\begin{theorem}
\label{thm2}
The following formula is valid:
\begin{equation}
\label{eq:alpharepr}
\chi'_G(3)= 3\ \sum ( - 1/3 )^{\rank(\FM(\alpha))/2}\, \sign ( \FM(\alpha));
\end{equation}
here the sum is calculated with respect to everywhere nonzero collections of spins $\alpha$ such that the rank of the matrix $\FM(\alpha)$ is even.
\end{theorem}

\subsection{Examples of calculations by formulas given in theo\-rems~\ref{thm1}~and~\ref{thm2}}
\label{subsec1-3}
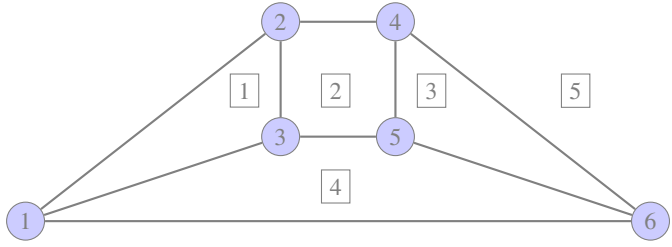
\begin{figure}[h]
\begin{center}
\begin{tikzpicture}
    \node[main node] (2) {$2$};
    \node[main node] (4) [right = 1cm  of 2]  {$4$};
    \node[main node] (3) [below = 1cm  of 2] {$3$};
    \node[main node] (5) [right = 1cm  of 3] {$5$};
    \node[main node] (1) [below left = 0.75cm and 3cm of 3]  {$1$};
    \node[main node] (6) [below right = 0.75cm and 3cm of 5]  {$6$};
    
    \node[draw]  [below left = 0.5cm and 0.1cm of 2]  {$1$};
    \node[draw]  [below right = 0.5cm and 0.35cm of 2]  {$2$};
    \node[draw]  [below right = 0.5cm and 0.1cm of 4]  {$3$};
    \node[draw]  [below right = 0.5cm and 2cm of 4]  {$5$};
    \node[draw]  [below right = 0.25cm and 0.35cm of 3]  {$4$};

    \path[draw,thick]
    (1) edge node {} (2)
    (1) edge node {} (3)
    (1) edge node {} (6)
    (2) edge node {} (3)
    (2) edge node {} (4)
    (3) edge node {} (5)
    (4) edge node {} (5)
    (4) edge node {} (6)
    (5) edge node {} (6)
        ;
\end{tikzpicture}
\caption{\label{pic:1} We will calculate the number of proper colorings of vertices and edges of this graph in 3 colors with the help of theorems~\ref{thm1}~and~\ref{thm2}, correspondingly.}
\end{center}
\end{figure}

Let us illustrate theorems~\ref{thm1}~and~\ref{thm2} by applying them for the calculation of the number of proper colorings of vertices in 3 colors (Example~1) and Tait colorings (Example~2) for the graph shown in Fig.~\ref{pic:1}. Recall that a coloring of vertices is said to be proper if adjacent vertices have different colors. 
The ratio between the number of such colorings for a connected graph and the number of colors~$q$ (in the case of a connected graph) equals the value of the characteristic polynomial of the cyclic matroid~$M_G$ at the point~$q$. 
The ground set $E$ of the cyclic matroid $M_G$ is the set of edges of the graph~$G$. 
The restriction $M_G|A$ represents a cyclic matroid of the subgraph induced by the set of edges~$A$;
this subgraph is not necessarily connected.
Recall that for a connected graph $G$ each spanning tree is a base of the matroid $M_G$. In the case of a nonconnected graph, the base of a cyclic set is the totality of spanning trees of all connectivity components. 

\textbf{Example~1.}
In order to apply Theorem~\ref{thm1}, let us make use (see, for exam\-ple,~\cite{stanleyArticle} and references therein) of the fact that the generating function of spanning trees $s(M_G;\alpha)$ is the algebraic complement of any element in the weighted Laplace--Kirchhoff matrix of the graph~$G$ (a generalized variant of the matrix theorem on trees). 
For the graph under consideration this matrix takes the form
$$
L=\left[
\begin{array}{cccccc}
 \alpha_{12}+\alpha_{13}+\alpha_{16} & -\alpha_{12} &
   -\alpha_{13} & 0 & 0 & -\alpha_{16} \\
 -\alpha_{12} & \alpha_{12}+\alpha_{23}+\alpha_{24} &
   -\alpha_{23} & -\alpha_{24} & 0 & 0 \\
 -\alpha_{13} & -\alpha_{23} &
   \alpha_{13}+\alpha_{23}+\alpha_{35} & 0 &
   -\alpha_{35} & 0 \\
 0 & -\alpha_{24} & 0 &
   \alpha_{24}+\alpha_{45}+\alpha_{45} & -\alpha_{45} &
   -\alpha_{46} \\
 0 & 0 & -\alpha_{35} & -\alpha_{45} &
   \alpha_{35}+\alpha_{45}+\alpha_{56} & -\alpha_{56} \\
 -\alpha_{16} & 0 & 0 & -\alpha_{46} & -\alpha_{56} &
   \alpha_{16}+\alpha_{46}+\alpha_{56} 
\end{array}
\right]  .
$$
This fact allows us to consider (with the help of a computer program) 512 values (with $\alpha_e\in \{1,-1\}\subseteq \mathbb F_3$) of the Legendre symbol $\left( \frac{\bar s(M_G;\alpha)}{3}\right)=\left( \frac{s(M_G;\alpha) \prod_e \alpha_e}{3}\right)$. In 164 cases, the algebraic complement equals zero. 
In 186 cases, $\left( \frac{\bar s(M_G;\alpha)}{3}\right)=+1$; in 162 cases, $\left( \frac{\bar s(M_G;\alpha)}{3}\right)=-1$.
In addition, the number of basic elements (edges of the spanning tree) equals~5, correspondingly, $r^*(M_G;\alpha)=|E|-5=4$, the contribution of $186+162$ variants in the right-hand side of formula~\eqref{eq:main1} equals $\frac{186}{9}-\frac{162}{9}$.

In order to consider the rest 164 variants, we need to know in these cases the value $s(M_G|A;\alpha)$ with $|A|=8$ (there exists, evidently, 9 such subsets~$A$). It coincides with the value $s(M_G;\alpha)|_{\alpha_e=0}$, where $\{e\}=E-A$. There are only 16 cases among 164 ones, when all these 9 values are equal to zero
(``zero cases'').
In the rest cases, $r^*(M_G;\alpha)=|A|-5=3$ and, consequently, these variants do not affect the right-hand side of formula~\eqref{eq:main1}.

In all 16 ``zero cases'' we have
$$
s_0(\alpha)\equiv s(M_G;\alpha)|_{\alpha_{12}=\alpha_{13}=0}\neq 0,\quad \text{i.e. for $A=E-\{(1,2),(1,3)\}$}\quad
\bar s(M_G|A;\alpha)=s_0(\alpha)\prod_{e\in A}\alpha_e\neq 0. 
$$ 
Here $r^*(M_G;\alpha)=2$. Applying the above formula for calculating $\bar s(M_G|A;\alpha)$, we conclude that $\left( \frac{\bar s(M_G|A;\alpha)}{3}\right)=1$ in 6 cases and
$\left( \frac{\bar s(M_G|A;\alpha)}{3}\right)=-1$ in 10 cases. As a result, we conclude that for the graph~$G$ shown in Fig.~\ref{pic:1},
$$
\chi_{M_G}(3)=\frac{186-162}{9}-\frac{6-10}{3}=4;
$$
consequently, there exist 12 proper colorings of vertices of this graph in 3 colors. Certainly, one can easily verify this fact by standard techniques developed for the calculation of chromatic polynomials.

\textbf{Example~2.}
Let us now calculate $\chi'_G(3)$ in accordance with Theorem~\ref{thm2}.
Let, for brevity, $(a,b,c,x,y,z)\equiv(\alpha(1),\alpha(2),\alpha(3),\alpha(4),\alpha(5) ,\alpha(6))$.
Then the matrix of faces for our graph takes the form
$$
\left[
\begin{array}{ccccc}
a+b+c & b+c     & 0      & a+c     & a+b  \\
b+c   & b+c+x+y & x+y    & c+y     & b+x  \\
0     & x+y     & x+y+z  & y+z     & x+z  \\
a+c   & c+y     & y+z    & a+c+y+z & a+z  \\
a+b   & b+x     & x+z    & a+z     & a+b+x+z
\end{array}
\right].
$$
All principal minors of the 4-th order coincide with each other and equal 
$
\sum_i \alpha(i_1)\alpha(i_2)\alpha(i_3)\alpha(i_4);
$
here the sum is calculated with respect to all quadruples of distinct spins, except for those with indices $\{2,3,4,5\}$, $\{1,3,5,6\}$ and $\{1,2,4,6\}$. Note that when deriving all formulas we have taken into account the fact that $\alpha(i)$ are nonzero elements of the field ${\mathbb F}_3$ and, consequently, $\alpha(i)^2=1$.
Calculating the sum in accordance with formula~\eqref{eq:alpharepr} with respect to 64 variants of values of $\alpha \in \{-1,1\}^{V}$, we conclude that the sum above in 30 cases gives a nonzero contribution; namely, in 18 cases this minor equals $-1 \pmod 3$, and in 12 cases it equals~1. In other 8 variants, this minor also equals zero, as well as all principal minors of the third order; moreover, several principal minors of the second order differ from zero, but in all these cases they equal $-1 \pmod 3$. 
Therefore, the sum in expression~\eqref{eq:alpharepr} equals
$$
\frac{-18+12}{9}-\frac{-8}{3}=2,
$$ 
i.e., $\chi'_G(3)=6$.

\section{Auxiliary assertions necessary for proving the\-o\-rems~\ref{thm1} and~\ref{thm2}}
\label{sect2}
\subsection{Remark on the set $A^*$}

In Theorem~\ref{thm1} we consider the set $A^*$ whose cardinal number is maximal among all subsets~$A$ of the set~$E$ such that $s(M|A)\neq 0$. We can essentially reduce the enumeration of all subsets. 
Here we prove one simple property of the set $A^*$, and later in Corollary~\ref{conc:1} reduce the enumeration even more essentially.

\begin{lemma}
\label{lem:rA*}
In denotations of Theorem~\ref{thm1} the following equality is valid:
$r(A^*)=r(M)$.
\end{lemma}

\begin{pf}
Assume that $s(M|A)\neq 0$ for a certain set~$A$ with the rank $r(A)<r(M)$.
Let the symbol $I$ stand for a base of the matroid~$M|A$, i.e. $|I|=r(A)$. According to the so-called independence augmentation property (see~\cite[subsection 1.1]{Oxley}),
there exists an element $e$ of the base of the matroid $M$ such that $r(I\cup e)=r(A)+1$. Consequently, $r(A\cup e)>r(A)$.

Let now $B$ be an arbitrary base of the matroid $M|(A\cup e)$. It necessarily contains $e$, otherwise for sets $A$ and $B$ we would get a contradiction with the second axiom of the rank function of the matroid (see subsection~\ref{subsec1-2}). Bur then 
$s(M|(A\cup e))=s(A)\times\alpha_e\neq 0$. Therefore the considered set $A$ does not coincide with $\mathop{\rm argmax} \{|A|: s(M|A)\neq 0 \}$,
i.e., $r(A^*)=r(M)$.
\end{pf}

\subsection{Another statement of Theorem~\ref{thm1}}

In paper~\cite{EJC}, we express the flow polynomial in terms of the linear combination of the Legendre symbols of values of the generating function of spanning trees of a connected graph. The flow polynomial is a characteristic polynomial of the bond matroid of the graph.
But spanning trees are not bases of the bond matroid, they represent bases of the cyclic matroid dual to the bond one. 
Note that in Theorem~\ref{thm1} the sum is calculated with respect to  bases of the matroid, whose characteristic polynomial is the desired one.
Here we state another assertion, whose matrix variant implies (as a particular case) the main result of the paper~\cite{EJC}.

Recall that
\textit{the dual matroid} $M^\perp=(E,r^\perp)$ with respect to the matroid $M=(E,r)$ is defined by the rank function $r^\perp$ which is defined 
on the same ground set as~$r$ and obeys the correlation
$$
r^\perp(A) = r(E - A) + |A| - r(E),\quad \mbox{for any $A\subseteq E$.}
$$
One can easily make sure that \textit{there is a biunique correspondence between bases of the initial matroid and the dual one}:
if $B$ is a base of $M$, then $E-B$ is a base of $M^\perp$ and vice versa. Hence $(M^\perp)^\perp=M$. It is also well known that 
(see, for example,~\cite{Oxley}) the matroid, which is dual to the matroid representable over the field $\mathbb F$, is also representable over this field. Consequently, the matroid dual to a regular one is also regular.
If $e\in E$ is a loop of the matroid $M^\perp$, then $e$ is said to be a {\it coloop} of the matroid $M$.

Let $A\subseteq E$, $C=E-A$.
\textit{The contraction} of $M$ onto $A$ is the matroid $M.A\equiv M/C$  with the ground set~$A$ and the rank function obeying the correlation
$$
r_{M.A}(B)= r_M (C \cup B) - r_M (C)\quad \mbox{for any $B\subseteq A$.}
$$

Assume that 
$$
A_\perp^*(\alpha)\equiv A_\perp^*=\mathop{\rm argmax}\limits_A \{|A|: s(M.A;\alpha)\neq 0 \}.
$$
The set $A_{\perp}^*$ is not defined uniquely, but the contribution of each term in sum~\eqref{eq:main2} is independent of the choice of this set.

Let us use the symbol $\mathbb F_q^*$ for the {\it set of nonzero elements} of the field $\mathbb F_q$ (they form a multiplicative group).
\begin{theorem}
\label{th:main2}
Let $q=p^d$ with odd prime $p$.
For any regular matroid without coloops the following formula is valid:
\begin{equation} \label{eq:main2}
\chi_{M^\perp} (q) =
\sum\nolimits_{\alpha\in (\mathbb F_q^*)^E}\ 
g(q,r(M.A_\perp^*(\alpha))\,
\eta(s(M.A_\perp^*;\alpha)).
\end{equation}
\end{theorem}

\begin{lemma}
\label{prop:1}
Theorem~\ref{th:main2} is equivalent to Theorem~\ref{thm1}.
\end{lemma}

\begin{pf}
As is well known, $(M|A)^\perp=M^\perp.A$, $M^\perp|A=M.A$. 
In addition, as was mentioned above, each base $B$ of the matroid $M|A$ corresp\-onds to the base
$(M|A)^\perp$, which equals $A-B$; therefore, $\bar s(M|A;\alpha)=s((M|A)^\perp;\alpha)$.
Thus, $A^*(M^\perp;\alpha)=A_\perp^*(\alpha)$.

Consequently, for establishing a biunique correspondence between each term in sums~\eqref{eq:main1} and~\eqref{eq:main2} (the matroid $M^\perp$ plays the role of the matroid $M$ in the first formula) it suffices to make sure that $r^*(M^\perp;\alpha)=r(M.A_\perp^*(\alpha))$.

Let us make use of Lemma~\ref{lem:rA*}. The following assertion is valid: $r_{M^\perp}(A^*_\perp)=r_{M^\perp}(E)$, i.e., $A^*_\perp$ contains a basis of the matroid~$M^\perp$.
Consequently, the comple\-ment of this set $E-A^*_\perp$ represents a part of the basis of the matroid~$M$, i.e.,
$r(E-A^*_\perp)=|E|-|A_\perp^*|$. 

We conclude that 
$$
r^*(M^\perp;\alpha)=|A_\perp^*|-r^\perp(E)=r(E)-|E|+|A_\perp^*|=r(E)-r(E-A_\perp^*)=r(M.A_\perp^*).
$$
\end{pf}

\subsection{A more general variant of Theorem~\ref{thm1}}

Theorem~\ref{thm1} is not valid for non-regular matroids. For example, the uniform matroid $U_{2,4}$, which, as is well known (see~\cite[Proposition 6.5.2]{Oxley}), is representable over all fields, except the field $\mathbb F_2$, satisfies the equality $\chi_{U_{2,4}}(q)=(q-1)(q-3)$.
At the same time, the right-hand side of formula~\eqref{eq:main1} equals 0 with $q=3$, but differs from 8 with $q=5$.

Note that matroid $U_{2,4}$ cannot be defined by a unimodular matrix 
over any field $\mathbb F_q$,
because it is not regular (see~\cite[Theorem 6.6.3]{Oxley}).
This matroid is representable by the matrix
\begin{equation}
\label{MU24}
\M(U_{2,4})=\left(\begin{matrix}1&0&1&1\\0&1&1&-1\end{matrix}\right).
\end{equation}
Evidently, the determinant of the submatrix formed by the last two columns in any field, whose characteristic exceeds 3, differs from $\pm 1$.

Let us now state a generalization of Theorem~\ref{th:main2}, which is valid for all matroids representable over the field~$\mathbb{F}_q$.
The statement of this generalization contains the representation matrix of the matroid.

Let $\M$ be the representation matrix for the matroid $M$ over the field $\mathbb F_q$; let the symbol~$E$ stand for the set of numbers of its columns, and let the symbol $V$ do for the set of numbers of its rows. 
Without loss of generality, we assume that the number of rows equals the rank of this matrix. 
Denote by $\M|_B$ the nondegenerate square submatrix of the matrix~$\M$ formed by columns whose numbers belong to the set~$B$
(and all rows of the matrix~$\M$). Put
$$
s'(\M;\alpha)=\sum_{B\in{\mathcal B}(\M)}{\det}^2 (\M|_B) \prod_{e \in B} \alpha_e,
$$
where ${\mathcal B}(\M)$ are all possible collections of columns of this matrix such that $\det (\M|_B)\neq 0$.

Let~$W\subseteq V$. Denote by $\M/W$ the matrix obtained from the matrix $\M$ by deleting rows whose numbers belong to the set~$W$; this matrix corresponds to some matroid. Evidently, the rank of this matroid equals $|V|-|W|$. 

Denote by $W^*({\M},\alpha)$ an arbitrary minimum cardinality subset $W$ of the set~$V$, for which the sum $s'(\M/W;\alpha)$ differs from zero. Denote by $r^*(\M;\alpha)$ the difference $|V|-|W({\M},\alpha)|$.

\begin{theorem}
\label{thm3}
Let us use the following denotations: $\mathbb F_q$ is a finite field of an odd characteristic, $M$ is an $\mathbb F_q$-linear matroid, and $\M$ is some matrix of its representation. The following formula is valid:
\begin{equation} \label{eq:main3}
\chi_{M^\perp} (q) =
\sum\nolimits_{\alpha\in (\mathbb F_q^*)^E}\ 
g(q,r^*(\M;\alpha))\,
\eta(s'(\M/W^*({\M},\alpha);\alpha)).
\end{equation}
In addition, (not only the sum, but also) each term in the sum depends neither on the choice of the representation matrix of the matroid~$\M$ nor on the choice of the collection~$W^*$ of rows of this matrix.
\end{theorem}

\begin{rmk}
In~\cite{EJC}, we consider a particular case of Theorem~\ref{thm3} for a cyclic matroid~$M$. We do not state the independence of each term in the sum \eqref{eq:main3} on the representation matrix of the matroid, but consider a fixed (canonic) matrix $\M$ that represents the incidence matrix of a connected graph (with the deleted row that corresponds to one of vertices). 
Since this matrix is unitary, we conclude that $\det^2 (\M|_B)\equiv 1$, and the sum calculated with respect to bases ${\mathcal B}(\M/W)$ in this case represents the sum with respect to spanning trees of the corresponding graph.
\end{rmk}

\textbf{Example 3.}
Let us illustrate correlation~\eqref{eq:main3} in the case of the matroid~$U_{2,4}$. Note that the matroid~$U_{2,4}$ is dual to itself, while the left-hand side of equality~\eqref{eq:main3} is
$\chi_{\,U_{2,4}^\perp}(q)=(q-1)(q-3)$. If~$U_{2,4}$ is defined by matrix~\eqref{MU24}, then
$$
s'(\M(U_{2,4});\alpha)=\alpha_1\alpha_2+\alpha_1\alpha_3+\alpha_1\alpha_4+\alpha_2\alpha_3+\alpha_2\alpha_4+4\alpha_3\alpha_4.
$$
We can prove that $r^*(\M;\alpha)=0$ only in $q-1$ cases (namely, when ${\alpha_1=\alpha_2}$, $\alpha_3=\alpha_4$, and $\alpha_1=-2\alpha_3$). Therefore, for~$U_{2,4}$ with representation matrix~\eqref{MU24} Theorem~\ref{thm3} is equivalent to the equality
$$
(q-1)(q-4)=g(q,2)\sum_{\alpha\in ({\mathbb F}_q^*)^4}  
\eta(\alpha_1\alpha_2+\alpha_1\alpha_3+\alpha_1\alpha_4+\alpha_2\alpha_3+\alpha_2\alpha_4+4\alpha_3\alpha_4).
$$

We see, that equality~\eqref{eq:main3} implies several new correlations for the Legendre symbol.

\subsection{Theorem~\ref{th:main2} as a corollary of Theorem~\ref{thm3} and some refinements with respect to the choice of $A^*$ and $A^*_\perp$}
\label{A*}
\begin{lemma}
\label{prop:2}
Theorem~\ref{thm3} implies Theorem~\ref{th:main2}.
\end{lemma}

\begin{pf}
Let $M$ be a regular matroid mentioned in assumptions of Theorem~\ref{th:main2}. According to Lemma~\ref{lem:rA*} and Proposition~\ref{prop:1}, $E-A^*_\perp(\alpha)$ is a part of a certain base of the matroid~$M$. Denote \textit{this base} by the symbol~$B^*$.  
Consider (for given fixed $\alpha$) the representation matrix $\M$ of the matroid~$M$, whose columns that correspond to the basis $B^*$ form the unit matrix. Note that due to the regularity of the matroid~$M$ we can choose the matrix $\M$ as a unimodular one. Therefore, the value of the function $s'$ for the matrix $\M$ equals the value of the function $s$ for the matroid that corresponds to this matrix.

Let $A'$ be some part of the base~$B^*$, which corresponds to columns, whose unit elements are located in rows, whose numbers belong to the set $W$. Then in accordance with~\cite[Proposition 3.2.6]{Oxley} by deleting rows $W$ and columns $A'$ from the matrix $\M$ we obtain the representation matrix of the matroid $M/A'$. If we delete only rows and do not delete columns, then we get a matrix with $|A'|$ zero columns, whose matroid differs from $M/A'$ in the existence of $|A'|$ additional loops. In both cases, the function $s$ is one and the same.

Therefore, for establishing a biunique correspondence between each term in sums~\eqref{eq:main2} and~\eqref{eq:main3} it suffices to make sure that in the unit submatrix of the matrix $\M$ composed of columns of~$B^*$ nonzero elements in rows~$W^*$ are located exactly in columns of $E-A^*_\perp$. The latter property is evident, because by construction $(E-A^*_\perp)\subseteq B^*$ and $A^*_\perp$ has the maximal cardinality among all such subsets~$A$ with nonzero $s(M.A)$. We choose the subset $W^*$ by using the same principle (naturally replacing the maximum with the minimum and doing $s$ with $s'$).
\end{pf}

Note that the number of all subsets $W$ of rows of the matrix, generally speaking, is less than the number of all possible subsets of the set~$E$; in the proof given above we fix some base and make use of the fact that each term in sum~\eqref{eq:main3} is independent of the representation matrix of the matroid. 

\begin{conclusion}
\label{conc:1}
If Theorem~\ref{thm3} is valid, then theorems~\ref{thm1} and~\ref{th:main2} remain valid, provided that the search of subsets $A^*$ and $A^*_\perp$ is performed by the enumeration of a lesser number of variants. Namely, let $B'$ be some fixed base of the matroid~$M$ mentioned in assumptions of theorems~\ref{thm1} and~\ref{th:main2}. Then one can define sets $A^*$ and $A^*_\perp$, correspondingly, by formulas 
$$
A^*=\mathop{\rm argmax}\nolimits_A \{|A|: B'\subseteq A,\ s(M|A;\alpha)\neq 0 \},
\qquad A_\perp^*=\mathop{\rm argmax}\nolimits_A \{|A|: (E-A) \subseteq B',\ s(M.A;\alpha)\neq 0 \}.
$$
\end{conclusion}

In particular, in Example~1 for 164 variants 
of values $\alpha$ we make sure that $s(M|A;\alpha)$ equals zero with each choice 
of $A\subseteq E$, $|A|=8$, among 9 possible ones. If we fix edges $B'$ of some spanning tree and consider only the set $A$ such that $B'\subseteq A$, then it suffices to verify the mentioned equality only for 4 values instead of 9 ones. Naturally, if we consider subsets $A$ of the set $E$ of a lesser cardinality, then the reduction of the enumeration process is even more essential.

\subsection{Flows and an analog of the Laplace--Kirchhoff matrix}
\label{subcec2-5}
In the paper~\cite{EJC}, when considering a particular case of Theorem~\ref{thm3}, we make use of properties of the Laplace--Kirchhoff matrix of the graph and the fact that the value of the flow polynomial at the point~$q$ by definition equals the number of everywhere nonzero solutions in the field $\mathbb F_q$ of the system of linear homogeneous equations with the incidence matrix of the graph. 
Let us now state analogs of these properties for an arbitrary matroid representable over the field~$\mathbb F_q$.
We use these properties in the proof of Theorem~\ref{thm3} in the next section.

Let us first recall the classical Capo--Rota result (\cite[Theorem~1 in sub\-sec\-tion~16.4]{crapo}) on the geometric sense of the value $\chi_M(q)$ for the matroid~$M$ such that columns of its representation matrix define $|E|$ points in the space~$\mathbb F_q^V$. As appeared, $\chi_M(q)$ equals the number of linear functionals in this space, which take on nonzero values at all these points.

J.~P.~S.~Kung \cite[Section 1]{kung3} understands {\it a flow} as a vector ${\mathbf f}=(f(e): e\in E, f(e)\in \mathbb F_q)$ such that $\M {\mathbf f}^T=0$. 
From the matrix $\M$ that defines the matroid~$M$ he has constructed another matrix (orthogonal to the matrix $\M$), which defines the dual matroid $M^\perp$. Then he has performed the bijection between everywhere nonzero flows for the initial matrix $\M$ and linear functionals, which take on nonzero values at points that define columns of the new matrix. Thus, J.~P.~S.~Kung has obtained the following assertion (\cite[Theorem~1.1, particular case]{kung3}, see also~\cite{kung4}) as a corollary of the Crapo--Rota result.

\begin{proposition}
\label{flow}
The number of everywhere nonzero flows for the matrix $\M$ equals $\chi_{M^\perp}(q)$.
\end{proposition}

Certainly, this assertion is well known in the case of a cyclic matroid~$M$. In this particular case, it is also well known that the principal minor of the weighted Laplace--Kirchhoff matrix, which can be expressed in terms of $\M$, represents a generating function of spanning trees $s(M;\alpha)$. 
Let us now recall the case of an arbitrary representable matroid.

Let $\M$ be the representation matrix of an arbitrarily chosen matroid~$M$. Recall that by our assumption the number of rows of this matrix equals $r(M)$. 
Consider the following square matrix
$L(\M;\alpha)$ of the order~$r(M)$:
\begin{equation}
\label{eq:LM}
L(\M;\alpha)=\M \Lambda \M^T,
\end{equation}
where $\Lambda$ is the $|E|\times|E|$-diagonal matrix whose diagonal elements equal $\alpha_e$,
$e\in E$, and $T$ is the transposition sign. One can easily prove the following assertion by using the Binet--Cauchy formula. 

\begin{proposition}[\cite{reiner1}, Exercise~11 (c)]
\label{prop:MatrixTheorem}
Let an $F_q$-linear matroid $M$ be defined by a matrix~$\M$. Then
$$
\det(L(\M;\alpha))=s'(\M,\alpha).
$$
\end{proposition}

\subsection{The Heawood theorem}
\label{subcec2-6}
For proving Theorem~\ref{thm2}, we (similarly to Proposition~\ref{flow}) need to represent the desired value (in this case, $\chi'_G(3)$) as the number of everywhere nonzero solutions to the system of linear equations over a finite field (in this case, over the field $\mathbb F_3$). This representation of the number of Tait colorings for a cubic graph was obtained by P.~J.~Heawood~\cite{heawood}. 

\begin{proposition}[Heawood, 1898]
\label{th:heawood}
Let us associate each vertex $v$ of an arbitrary planar biconnected cubic graph~$G$ with the variable (spin) $\sigma(v)$, which takes on values in the set $\{-1,1\}\equiv F_3^*$. 
Then the normed number of Tait colorings $\chi'_G(3)/3$ equals the number of all possible sets of spins $(\sigma(v), v\in V(G))$ such that for any face~$F$ of this graph the sum of spins for all vertices of the face~$F$ equals zero.
\end{proposition}

In the original paper by P.~J.~Heawood, this assertion is stated in terms of the graph dual to~$G$ (see also \cite{belaga} for another proof of the same theorem). However, it is often more convenient to consider a cubic graph (cf. \cite[Theorem 9.3.4]{Ore}).

\subsection{Multidimensional Gaussian sums over the field ${\mathbb F}_q$}
\label{subcec2-7}
In the application of the $\alpha$-representation in the case of a real field, explicit formulas for calculating Gaussian integrals with the imaginary unit in the expo\-nent play an important role. In the case of a finite field, we mean multidimensional Gaussian sums. Here we also need these formulas (established in \cite{EJC}). 

The function $\exp (i x)$ is a homomorphism of an additive group of real numbers to the group of complex numbers which are modulo equal to one; all the rest such homomorphisms (additive characters) are parameterized by the parameter $k\in R$ and take the form $\exp (i k x)$.
In the case of a finite field $\mathbb F_q$, $q=p^d$, an analogous role is played by the function 
\begin{equation}
\label{eq:h}
h(x)=\exp{(2\pi i\, \operatorname{Tr}(x)/p)}, \mbox{ where $\operatorname{Tr}(x)=x+x^p+x^{p^2}\ldots+x^{p^{d-1}}$;}
\end{equation}
all additive characters take the form $h_k(x)=h(k x)$, $k\in \mathbb F_q$ (\cite[Theorem~5.7]{lidlNider}).

Let $C$ be an arbitrary symmetric $n\times n$ matrix  with elements belonging to~$\mathbb F_q$, where $q$ is odd, while ${\mathbf x} C {\mathbf x}^T$ is a quadratic form with this matrix (${\mathbf x}$ is a row vector of the corresponding dimension). 
The following sum is an analog of the multidimensional Gaussian integral: 
$$
\Gau_q(C)=\sum\nolimits_{{\mathbf x}\in \mathbb F_q^n} h({\mathbf x} C {\mathbf x}^T).
$$
In the one-dimensional case with $C=1$ we consider the quadratic Gaussian sum $g_1(q)=\sum_{x\in\mathbb F_q} h(x^2)$. As is well known (\cite[Theorem~5.15]{lidlNider}), for a field with an odd characteristics it obeys the formula
$$
g_1(q) = \left\{ \begin{array}{ll}
   (-1)^{d-1}\sqrt{q},& \textrm{if}~p \bmod 4= 1,  \\
   (-1)^{d-1} i^d \sqrt{q},&  \textrm{if}~p \bmod 4=3.  \end{array} \right.
$$
Note that the function $g(q,m)$ defined by us earlier with even $m$ coincides with $\left[ {\frac{g_1(q)}{q}} \right]^m$.
Below we give an explicit formula for $\Gau_q(B)$, but first let us make one useful remark.

\begin{rmk}
\label{rem:cong}
If matrices $C$ and $A$ are congruent, i.e., $A=P C P^T$, where $P$ is a nondegenerate $n\times n$ matrix, then $\Gau_q(C)=\Gau_q(A)$.
\end{rmk}

Remark~\ref{rem:cong} is valid, because with ${\mathbf x}'= P {\mathbf x}$ the sum $\Gau_q(A)$ turns into the sum $\Gau_q(C)$.

\begin{lemma}[Lemma~8 in~\cite{EJC}]
\label{lem:gauss}
Let $q=p^d$ with odd prime $p$, $\operatorname{rank}\, C = r$. Then
\begin{equation}
\label{eq:gauss}
\frac{\Gau_q(C)}{q^n}  =
 \eta(\det C_r) \left[ {\frac{g_1(q)}{q}} \right]^r,
\end{equation}
where $\det C_r$ is an arbitrary nonzero principal minor of the order~$r$.
\end{lemma}

In view of Remark~\ref{rem:cong} the proof of Lemma~\ref{lem:gauss} is reduced to the diagonal case (see~\cite[Chapters IV]{serre} for the technique for reducing a quadratic form over a finite field to the diagonal form). 
The diagonal case is factorized to the one-dimensional variant, which in turn is reduced to the calculation of the sum $g_1(q)$.

\begin{conclusion}
\label{concl2}
For any symmetric matrix of the rank~$r$ the value $\eta(\det C_r)$, where $\det C_r$ is a nonzero principal minor of the matrix $C$ of the order~$r$, is independent of $C_r$.
\end{conclusion}

\begin{rmk}
In the case of zero matrix~$C$ the formula obtained in Lemma~\ref{lem:gauss} is valid, provided that in this case $\eta(\det C_r)=1$.
\end{rmk}

\begin{conclusion}
\label{concl3}
Assume that all elements of a symmetric matrix $C$ represent linear functions with respect to some collection of variables $\alpha\in\mathbb (F_q^*)^k$, while $r(C(\alpha))$ is the rank of this matrix. 
Then 
$
\sum_{\substack{
\alpha:\ \alpha\in (F_q^*)^k, \\
r(C(\alpha))\bmod 2 = 1}} \Gau_q(C(\alpha)) =0.
$
\end{conclusion}

\begin{pf}
Let $\gamma$ be an arbitrary element of the field~$\mathbb F_q$ such that $\eta(\gamma)=-1$.
Let us replace the collection $\alpha$ in the sum under consideration with $\gamma\alpha$. Note that then in the formula for $\Gau_q(C(\gamma \alpha))$ we get the value $\gamma^r \det C_r$ in place of $\det C_r$. For odd $r$ we get the equality $\eta(\gamma^r \det C_r)=-\eta(\det C_r)$.
Consequently, the considered sum equals itself with the opposite sign.
\end{pf}

\section{The $\alpha$-representation}
\label{sec3}
\subsection{The Fourier transform over the field ${\mathbb F}_q$ and its properties}
\label{bog}
For proving theorems~\ref{thm2} and~\ref{thm3} we use several simple properties of the Fourier transform over the field~${\mathbb F}_q$. 
Consider complex-valued functions $f(k)$ whose argument $k$ belongs to the field $\mathbb{F}_q$.
We understand the Fourier transform of such a function as the function $\widehat f(x)=\sum_{k\in \mathbb{F}_q} f(k)h(kx)/q$, where $h$ obeys formula~\eqref{eq:h}.

Note that usually the definition of the Fourier transform proposed by us is used for the inverse discrete Fourier transform. 
It is more convenient for us to use the terminology proposed in the paper~\cite{EJC}. The use of the classical definition would require to introduce the conjugation sign and an additional factor $q$ in the right-hand sides of the corresponding formulas. 

Let the symbol ${\mathbf 1}(k)$ denote the function $f(k)$ that identically equals one, and let the symbol $\delta(x)$ stand for the delta function (the Kronecker symbol), i.e.,
$\delta(0)=1$, $\delta(x)=0$ with all $x\in\mathbb F_q^*$. 
It is well known that $\sum_{k\in \mathbb{F}_q} h(kx)=0$ for any nonzero $x$ (\cite[Theorem~5.4]{lidlNider}).
Hence
\begin{equation}
\label{widehat1}
\widehat{{\mathbf 1}}(x)=\delta(x). 
\end{equation}
Correlation~\eqref{widehat1} implies the formula
\begin{equation}
\label{alpha}
\sum\nolimits_{x\in \mathbb{F}_q^*} h(kx) =q \delta(k)-1=q\delta(k^2)-1=\sum\nolimits_{y\in \mathbb{F}_q^*} h( k^2 y);
\end{equation}
we use it in what follows.

\begin{rmk}
\label{rem:alpha}
In the case of a finite field, we can treat the function $f(k)=1-\delta(k)$ as the norm of an element of a finite field and do the sum $\sum_{x\in \mathbb{F}_q^*} h(kx)$ as the Fourier transform of the norm raised to a certain power.
In the case of a real field, an analog of this sum is the Fourier transform of (a generalized function) $|k|^\gamma$.
In the quantum field theory, such functions are known as the so-called propagators of Feynman amplitudes (in the massless case).
The parametric representation for integrals of propagators of Feynman amplitudes as integrals of characters with a quadratic argument was proposed by R.~Feynman. The paper~\cite[p.~691]{syman} by K.~Symanzik has given raise to their systematic application; the symbol $\alpha$ stands there for an analog of the variable~$y$ used by us.
In the quantum field theory, this representation of a propagator is called the $\alpha$-re\-pre\-sen\-ta\-ti\-on.
Following this tradition, we recall that Theorem~4 is a direct generalization of the main result obtained in the paper~\cite{EJC}; 
it has a real analog in the theory of Feynman amplitudes.
\end{rmk}

\subsection{Proof of Theorem~\ref{thm2}}

For clarifying the idea of using the $\alpha$-representation (in a discrete case, this idea is rather evident), it makes sense to consider a more concrete case first.

\begin{pot2}
According to the Heawood theorem (Proposition~\ref{th:heawood}), $\chi'_G(3)$ equals the tripled value of the sum
\begin{equation}
\label{eq:S}
S=\sum_{\sigma\in\{-1,1\}^{V(G)}} \prod_{F\in \mathcal F} \delta(\sum_{v\in F} \sigma(v)),
\end{equation}
where $\mathcal F$ is the {\it set of all faces of the graph~$G$}. 

Let us transform the right-hand side of formula~\eqref{eq:S}, using the fact that each $\delta$-function represents the Fourier transform of the unit (see~\eqref{widehat1}). 
Representing the product of exponents as the exponent of a sum and changing the summation order, we obtain the correlation
$$
S=\sum_{{\mathbf k}\in {\mathbb F}_3^{\mathcal F}} \sum_{\sigma\in\{-1,1\}^{V(G)}} \exp\left(\frac{2\pi i}3 \, \sum_{F\in\mathcal F} k_F\sum_{v\in F} \sigma(v) \right)/3^{|\mathcal F|}.
$$
We can represent the sum in the exponent in another way, namely, 
$$
\sum_{F\in\mathcal F} k_F\sum_{v\in F} \sigma(v) =\sum_{v\in V(G)} \sigma(v) \sum_{F: v\in F} k_F.
$$
This allows us to apply formula~\eqref{alpha}. We obtain the correlation
$$
S=\sum_{{\mathbf k}\in {\mathbb F}_3^{\mathcal F}} \sum_{\alpha\in\{-1,1\}^{V(G)}} \exp\left(\frac{2\pi i}3 \, \sum_{v\in V(G)} \alpha(v) \left( \sum_{F: v\in F} k_F\right)^2\ 
\right)/3^{|\mathcal F|}.
$$
By changing the order of outer sums and combining similar terms in the exponent, we obtain the following correlation:
$$
S=\sum_{\alpha\in\{-1,1\}^{V(G)}} \sum_{\mathbf k\in {\mathbb F}_3^{\mathcal F}} \frac{ \exp( 2\pi i\ {\mathbf k} \FM(\alpha)\,{\mathbf k}^T\ /3)}{3^{|\mathcal F|}}
=\sum_{\alpha\in\{-1,1\}^{V(G)}} \frac{\Gau_3(\FM(\alpha))}{3^{|\mathcal F|}}
$$
(here $\mathbf k$ is the vector row $(k_F, F\in\mathcal F)$).
Let us now make use of Corollary~\ref{concl3}. The right-hand side of formula~\eqref{eq:gauss} with $q=3$ and an even rank of the argument the function $\Gau_3$ coincides 
with the term in sum~\eqref{eq:alpharepr}.
\end{pot2}  

\subsection{Proof of Theorem~\ref{thm3}}

The proof of Theorem~\ref{thm3} proposed by us in many respects is analogous to the proof of Theorem~\ref{thm2}. Note also the similarity of final results.

In accordance with formula~\eqref{eq:LM} and Proposition~\ref{prop:MatrixTheorem} the value $s'(\M/W;\alpha)$ represents the determinant of the submatrix of the matrix $L(\M;\alpha)$ obtained from the latter by deleting rows and columns, whose numbers belong to the set~$W$. Correspondingly, by Lemma~\ref{lem:gauss} each term in sum~\eqref{eq:main3} multiplied by $q^{r(M)}$ equals (with even rank of the matrix $L(\M;\alpha)$) the Gaussian sum $\Gau_q(L(\M;\alpha))$. 

Recall that numbers of rows of the matrix $\M$ belong to a certain set~$V$, $r(M)=|V|$.

\begin{pot3}
Proposition~\ref{prop:1} implies the formula
$$
\chi_{M^\perp}(q)=\sum_{{\mathbf x}\in ({\mathbb F}_q^*)^E} \prod_{j\in V} \delta(\M_j\, {\mathbf x}^T),
$$
where $\M_j$ is the $j$-th row of the matrix $\M$, ${\mathbf x}=(x_e, e\in E)$. 

Similarly to the proof of Theorem~\ref{thm2} we use formula~\eqref{widehat1} and, transforming the product of additive characters~$h$ and changing the order of summation, obtain the correlation
$$
\chi_{M^\perp}(q)=\sum_{{\mathbf k}\in {\mathbb F}_q^V} \sum_{{\mathbf x}\in ({\mathbb F}_q^*)^E} h\left(\sum_{j\in V} k_j \M_j\, {\mathbf x}^T \right)/q^{|V|}.
$$
Let us write down the sum in the argument of the function~$h$ termwisely by $x_e$, namely,
$$
\sum_{j\in V} k_j \M_j\, {\mathbf x}^T =\sum_{e\in E} x_e \sum_{j\in V} k_j \M_{j,e}
$$
This allows us to apply formula~\eqref{alpha} and thus obtain the correlation
$$
\chi_{M^\perp}(q)=\sum_{{\mathbf k}\in {\mathbb F}_q^V} \sum_{\alpha\in ({\mathbb F}_q^*)^E} 
h\left(\sum_{e\in E} \alpha_e \left( \sum_{j\in V} k_j \M_{j,e}\right)^2 \right)/q^{|V|}.
$$
Let us again change the order of outer sums and transform the argument of the function~$h$. As a result, we conclude that
$$
\chi_{M^\perp}(q)=
\sum_{\alpha\in ({\mathbb F}_q^*)^E} \sum_{{\mathbf k}\in {\mathbb F}_q^V} \frac{ h( {\mathbf k}\, L(\M;\alpha)\,{\mathbf k}^T)}{q^{|V|}}
={\sum_{\alpha\in ({\mathbb F}_q^*)^E} \Gau_q(L(\M;\alpha))}/{q^{|V|}}.
$$
In accordance with Corollary~\ref{concl3} we can neglect terms with an odd rank of the matrix $L(\M;\alpha)$.

One can easily prove that each term in the latter sum is independent of the choice of the representation matrix of the matroid~$M$. Really, let $\M'$ be some other representation matrix of the same matroid. Then $\M'=P \M$, where $P$ is some nondegenerate matrix.
But then the matrices $L(M;\alpha)=\M \Lambda \M^T$ and $L(M';\alpha)$ are congruent, consequently, according to Remark~\ref{rem:cong}, 
$$\Gau_q(L(\M;\alpha))=\Gau_q(L(\M';\alpha)).$$
\end{pot3}

\section{Conclusion} 
This paper is devoted to the extention of the applicability domain of the $\alpha$-representation technique, which was first used in~\cite{EJC} for combinatorial problems, to the case of general structures of the enumerative combinatorics. To this end, we first express the number of everywhere nonzero solutions of a system of homogeneous equations over a finite field $\mathbb F_q$ as the sum of analogs of multidimen\-si\-o\-nal Gaussian integrals over~$\mathbb F_q$. The latter allow an explicit representation in terms of the Legendre symbol of principal minors of the matrix in the quadratic form. As a result, we obtain a linear combination of these Legendre symbols with rational coefficients which modulo equal natural exponents of $1/q$. 

We have considered several variants of these expressions for the characteristic polynomial of the matroid. 
Note that in preprint~\cite{arXiv} we propose another (more complicated) proof of the simplified variant of Theorem~\ref{thm3}; there we use results obtained by C.~Chevalley for the number of values of variables in the field $\mathbb F_q$, with which the quadratic form equals zero.

We have also obtained the $\alpha$-representation for the number of Tait colorings for planar cubic graphs.
Note that in Example~2 the equality of all principal minors of the 4-th order is due to standard properties of the Laplace--Kirchhoff matrix. In our further paper, we are going to represent the matrix of faces as a weighted Laplace--Kirchhoff matrix for the dual graph and to explicitly describe variants of weights of its edges, which are taken into account in the sum. 
Certainly, this research was inspired by its connection with the Four Colors Theorem.

\subsection*{Acknowledgements} 
The author is grateful to all participants of the grant online workshop for useful remarks on this work.
The author is also grateful to Olga Kashina for her help in the translation of this paper.

\end{document}